\documentclass[smallcondensed]{svjour3}
\usepackage{amsfonts}
\usepackage{latexsym}
\usepackage{amsmath}
\usepackage{amssymb}
\usepackage{amscd}
\usepackage{epsfig}
\usepackage[affil-it]{authblk}
\addtocounter{MaxMatrixCols}{4}

\begin{document}

\title{A Note on Using the Resistance-Distance Matrix to solve Hamiltonian Cycle Problem}
\author{V. Ejov, J.A. Filar, M. Haythorpe$^*$, J.F. Roddick, and S. Rossomakhine}
\institute{V. Ejov
\at Flinders University\\
1284 South Road, Clovelly Park SA 5042 \\
\email{vladimir.ejov@flinders.edu.au}
\and
J.A. Filar
\at Flinders University\\
1284 South Road, Clovelly Park SA 5042 \\
\email{jerzy.filar@flinders.edu.au}
\and
M. Haythorpe - {\em Corresponding Author}
\at Flinders University\\
1284 South Road, Clovelly Park SA 5042 \\
\email{michael.haythorpe@flinders.edu.au}
\and
J.F. Roddick
\at Flinders University\\
1284 South Road, Clovelly Park SA 5042 \\
\email{john.roddick@flinders.edu.au}
\and
S. Rossomakhine
\at Flinders University\\
1284 South Road, Clovelly Park SA 5042 \\
\email{serguei.rossomakhine@flinders.edu.au}
}

\maketitle {\abstract An instance of Hamiltonian cycle problem can be solved by converting it to an instance of Travelling salesman problem, assigning any choice of weights to edges of the underlying graph. In this note we demonstrate that, for difficult instances, choosing the edge weights to be the resistance distance between its two incident vertices is often a good choice. We also demonstrate that arguably stronger performance arises from using the inverse of the resistance distance. Examples are provided demonstrating benefits gained from these choices.}

\section*{Acknowledgements}

This work was supported by Australian Research Council grants LP110100166 and DP150100618.

\section{Introduction}\label{sec-Introduction}

The famous {\em travelling salesman problem} (TSP) can be summarised as follows: given a graph $G$ containing a vertex set $V$ and an edge set $E : V \rightarrow V$, and a weighting function $f$ that assigns a weight to each edge in the graph, find a simple cycle that includes every vertex in $V$ such that the sum of weights of the edges used is minimised. TSP is known to be NP-hard, and the decision variant, in which the question is whether any tour with length less than a given value $k$ exists, is NP-complete. In the language of TSP, a simple cycle visiting every vertex is called a {\em tour}, however it is also known as a {\em Hamiltonian cycle}. A closely related problem is the {\em Hamiltonian cycle problem} (HCP) which simply asks if at least one Hamiltonian cycle exists in a given graph. Graphs with at least one Hamiltonian cycle are called {\em Hamiltonian} and graphs with none are called {\em non-Hamiltonian}.

Although specialised algorithms and heuristics for solving HCP do exist (e.g. see Baniasadi et al. (2014); Chalaturnyk (2008); Eppstein (2003)), a common and typically quite successful method for solving HCP is to first cast it as a TSP problem, and then use one of the suite of highly-developed TSP heuristics (e.g. see Applegate et al. (2006); Helsgaun (2000)) to solve it. It is natural to convert HCP to TSP, since in the latter a tour of minimal length is desired, and any tour suffices to give an affirmative answer to an instance of HCP. Specifically, HCP is converted to an optimisation (typically, integer programming) problem of the generic form:

\begin{eqnarray}&\min\limits_{\bf x} \sum\limits_i \sum\limits_j w_{ij} x_{ij}, \;\;\; \mbox{s.t. } \;\; {\bf x} = \{x_{ij} \in \{0, 1\} \;\;|\;\; (i,j) \in E\}& \label{eq-1}\end{eqnarray}

where the positive solution values corresponds to the edges in a tour, and $w_{ij}$ is the weight or cost associated with edge $(i,j)$. However, there is the question of how best to assign the weights.

One approach which, perhaps, seems logical is to assign every edge in the graph an equal weight, say, 1. This avoids unnaturally biasing the heuristic used by ensuring that all tours have equal weight, equal to $|V|$. However, this naive approach can cause the algorithm to run inefficiently. Since all edges are weighted equally, it is difficult for any TSP heuristic to make good choices, and often random guesses are required to find a solution. Using randomly chosen weights for the edges can help drive the heuristic in some direction, and if the graph has many tours this may be sufficient, but if there are relatively few tours this often amounts to little better than a prohibitively exhaustive search.

Since randomly generated graphs typically have exponentially many tours (and in fact, can be solved in general in only slightly longer than linear time (Frieze 2015)), graphs with relatively few tours usually contain significant structure. In this note, we propose to use information about the structure present in the underlying HCP instance to provide a better choice of weighting function. Specifically, we use information obtained from the {\em resistance-distance matrix} that contains all of the resistance distances between pairs vertices in a graph. Using the state-of-the-art TSP solver {\em Concorde} (Applegate et al. 2006), we attempt to solve several difficult families of HCP instances, and demonstrate that these instances can be much more easily solved using resistance distances. Interestingly, we show that for many families of HCP instances, the best choice is to use the inverse of resistance distances.

\section{Resistance-Distance Matrix}

Suppose that a graph $G$ corresponds to an electrical network, such that the vertices on the graph are electrical components, and each edge $(i,j)$ corresponds to a 1 ohm resistor being placed between components $i$ and $j$. The resistance distance $\Omega_{ij}$ is then the effective resistance between components $i$ and $j$, and the matrix of all resistance distances is called the resistance-distance matrix.

Although there are several ways to formulate resistance distance mathematically, perhaps the simplest arises from the inverse of the {\em Laplacian matrix}. Suppose that a graph $G$ of order $n$ has a corresponding adjacency matrix $A$, and further define a diagonal matrix $D$ where each entry $d_{ii}$ is equal to the degree of vertex $i$ in $G$. Then the Laplacian matrix $L := D - A$. Define a new matrix $\Gamma := L + \frac{1}{n} J$, where $J$ is the $n \times n$ matrix with every entry equal to 1. Then the resistance distance between vertices $i$ and $j$ was shown in (Babic et al. 2002) to be equal to

$$\Omega_{ij} = \left(\Gamma\right)^{-1}_{ii} + \left(\Gamma\right)^{-1}_{jj} - 2\left(\Gamma\right)^{-1}_{ij}.$$

Since an inverse needs to be found, it takes $O(n^3)$ time to compute the full set of resistance distances for every pair of vertices in the graph. However, for difficult instances of HCP this is a relatively small expense compared to the time taken to solve the instance. Note that the above definition is only well-defined if $G$ is connected. If $G$ is disconnected then the resistance distance can be found between vertices in connected components by treating those connected components as individual graphs, and the resistance distance between vertices in different connected components is defined to be infinity (see p.132 in Bapat (2010)).

Research into the value of the resistance-distance matrix is still relatively young, arguably initiated in 1993 by Klein and Randi\'{c} (1993), although the underlying theory of resistive electrical networks has been studied for far longer (e.g. see Doyle and Snell (1984)). The resistance-distance matrix has been used extensively in mathematical chemistry, where it is also known as the commute-time distance. For example, it is used as a tool in studying cyclicity and stability of molecules (e.g. see Babi\'{c} et al. (2002); Fowler (2002)) as well as signal transduction in proteins (e.g. see Chennubhotla and Bahar (2007)). Arguably the most famous result arising from this line of research is the Kirchhoff index (Lukovits et al. (1999)), defined as one half of the sum of all entries of the resistance-distance matrix. Palacios (2001) computed closed-form formulae for the Kirchhoff index of various classes of graphs. In the years since it was introduced, the resistance-distance matrix has been recognised as a useful tool in analysing graphs and networks of various kinds, including communications networks (Tizghadam et al. 2010), social networks (Alguliev et al. 2011), ecological connectivity models (McRae et al. 2007), and graphs arising from natural language processing (Rao et al. 2008), while Spielman and Srivastava (2011) use a nearly linear-time approximation of the resistance-distance matrix to develop a fast graph sparsification algorithm.

We now extend the use of the resistance distance matrix to aid in solving a classical graph theory problem, namely, HCP. In particular, in the formulation (\ref{eq-1}), we consider $w_{ij} = \Omega_{ij}$ and $w_{ij} = \Omega_{ij}^{-1}$, as well as randomly generated $w_{ij}$'s. In the following section, we demonstrate that difficult classes of HCP instances can be much more effectively solved by using resistance distances to weight the edges. This finding is, perhaps, unsurprising since the resistance distance is in some sense a measure of how many paths exist between two vertices. Intuitively then, it would seem to be linked to the problem of trying to find a longest cycle. Of course, such a link need not be direct, and so we will consider both the resistance distance, and the inverse of the resistance distance, and will show that in various cases, one or both are more effective than random weights.

\section{Examples}

We now provide examples of families of HCP instances which, in our experience, are difficult to solve, primarily due to the relative rarity of Hamiltonian cycles.

\begin{itemize}
\item {\bf Aldred-Thomassen graphs} - Aldred and Thomassen (1999) produced a family of 3-connected graphs which has only a single Hamiltonian cycle, and for which all vertices are degree 3 except two vertices of degree 4. One example of such a graph is displayed in Aldred and Holton (1999) and can be generalised easily. The graphs have order $16 + 4k$ for integer $k \geq 1$.
\item {\bf Fleischner-2 and Fleischner-3 graphs} - Fleischner (2014) provides two infinite families of graphs which all have minimum degree 4, and only a single Hamiltonian cycle. The first family contains 2-connected graphs of order $169k$ for integer $k \geq 2$. The second family contains 3-connected graphs of order $85 + 323k$ for integer $k \geq 2$. We will henceforth refer to the two families as Fleischner-2 and Fleischner-3, respectively.
\item {\bf Modified Flower Snarks} - Isaacs (1975) constructed an infinite family of Snarks, that is, non-Hamiltonian 3--regular graphs with chromatic index 4, called the Flower snarks. The family contains graphs of order $8 + 4k$ for any odd integer $k \geq 1$, which are constructed by taking $k$ copies of the Star graph and joining them together in a prescribed way. The graphs are all {\em maximally non-Hamiltonian}, that is, the addition of any edge renders them Hamiltonian. We modify the Flower snarks by adding a single edge. Specifically, we add an edge between any two non-adjacent vertices in one of the star graphs (any such choice is equivalent).
\item {\bf Minimally-Hamiltonian Regular Graphs} - Haythorpe (2016) proposed a family of regular graphs which are conjectured to contain the minimal number of Hamiltonian cycles over all $k$-regular graphs (for $k \geq 5$) of equivalent order. The family contains graphs of order $m(k+1)$ for integer $m \geq 2$ and $k \geq 5$. The graphs contain $(k-1)^2\left[(k-2)!\right]^{m}$ Hamiltonian cycles. In our experiments we chose $k = 10$ in all cases.
\end{itemize}

We considered several instances of graphs of the above types. In each case, we attempted to solve the graph by submitting it to {\em Concorde} in sparse format with edge weights chosen from the following four schemes (rounding up when necessary):

\begin{enumerate}\item All edge weights are equal to 1.
\item All edge weights are randomly chosen integers between 1 and 100.
\item All edge weights are equal to the resistance distance between their two incident vertices multipled by 100.
\item All edge weights are equal to the inverse of the resistance distance between their two incident vertices multiplied by 100.\end{enumerate}

In Table \ref{tab-results} we demonstrate how quickly Concorde was able to solve the various instances for the four schemes. In each case, we limited the size of the instance to less than 2000 vertices to ensure memory management was not an overriding factor in the experiment. We stopped Concorde if an individual run took longer than 24 hours. The experiment was conducted on a Linux machine with an AMD Opteron 6282 SE 2.6GHz Processor and 512GB RAM.

\begin{table}[h!]\begin{center}\begin{tabular}{|l|c||c|c|c|c|} \hline
{ {\bf Family} } & {  {\bf Order} } & { {\bf Scheme 1}} & {  {\bf Scheme 2}} & {  {\bf Scheme 3}} & {  {\bf Scheme 4}}\\
\hline { Aldred-Thom. } & {  736 } & {  00:07:52 } & {  00:02:21 } & {  00:02:35 } & {  {\bf 00:00:01}} \\
\hline { Aldred-Thom. } & {  976 } & {  00:07:19 } & {  00:03:20 } & {  00:03:32 } & {  {\bf 00:00:01}} \\
\hline { Aldred-Thom. } & {  1216 } & {  Failed } & {  00:04:11 } & {  00:10:19 } & {  {\bf 00:00:02}} \\
\hline { Aldred-Thom. } & {  1536 } & {  Failed } & {  00:06:37 } & {  00:06:29 } & {  {\bf 00:00:02}} \\
\hline \hline { Fleischner-2 } & {  338 } & {  Timeout } & {  00:06:08 } & {  00:02:29 } & {  {\bf 00:00:01}} \\
\hline { Fleischner-2 } & {  507 } & {  Timeout } & {  00:37:52 } & {  00:06:28 } & {  {\bf 00:00:01}} \\
\hline { Fleischner-2 } & {  676 } & {  Timeout } & {  Timeout } & {  00:21:42 } & {  {\bf 00:00:01}} \\
\hline { Fleischner-2 } & {  845 } & {  Timeout } & {  Timeout } & {  01:08:31 } & {  {\bf 00:00:01}} \\
\hline \hline { Fleischner-3 } & {  731 } & {  Timeout } & {  Timeout } & {  00:19:08 } & {  {\bf 00:00:01}} \\
\hline { Fleischner-3 } & {  1054 } & {  Timeout } & {  Timeout } & {  Timeout } & {  {\bf 00:00:02}} \\
\hline { Fleischner-3 } & {  1377 } & {  Timeout } & {  Timeout } & {  Timeout } & {  {\bf 00:00:02}} \\
\hline { Fleischner-3 } & {  1700 } & {  Timeout } & {  Timeout } & {  Timeout } & {  {\bf 00:00:02}} \\
\hline \hline { Flower } & {  324 } & {  Timeout } & {  Timeout } & {  {\bf 00:00:01} } & {  Timeout} \\
\hline { Flower } & {  684 } & {  Timeout } & {  Timeout } & {  {\bf 00:00:02} } & {  Timeout} \\
\hline { Flower } & {  996 } & {  Timeout } & {  Timeout } & {  {\bf 00:00:05} } & {  Timeout} \\
\hline { Flower } & {  1332 } & {  Timeout } & {  Timeout } & {  {\bf 00:00:06} } & {  Timeout} \\
\hline \hline { Minimal Reg. } & {  1100 } & {  00:08:52 } & {  00:00:25 } & {  00:06:44 } & {  {\bf 00:00:14}} \\
\hline { Minimal Reg. } & {  1375 } & {  04:34:59 } & {  00:00:30 } & {  00:10:40 } & {  {\bf 00:00:21}} \\
\hline { Minimal Reg. } & {  1650 } & {  09:07:27 } & {  00:00:58 } & {  00:13:26 } & {  {\bf 00:00:15}} \\
\hline { Minimal Reg. } & {  1925 } & {  21:09:23 } & {  00:02:06 } & {  00:27:10 } & {  {\bf 00:00:26}} \\
\hline \end{tabular}\end{center}\caption{Times taken (hh:mm:ss) to solve each instance in Concorde using the three weighting schemes given above. The fastest runtime for each instance is given in bold. Any run that took longer than 24 hours was terminated (labelled \lq\lq Timeout"). The label \lq\lq Failed" means that Concorde encountered an unexpected bug and crashed during execution.\label{tab-results}}\end{table}

It is interesting to note the vast range of solving times, depending on the scheme used. For example, the Fleischner graphs with over 1000 vertices all took longer than 24 hours to solve with three of the schemes, but mere seconds with Scheme 4. In fact, Scheme 4 typically dominated the other methods, but was unable to solve the modified Flower Snarks which in turn were solved in seconds by Scheme 3. These findings highlight the importance of using appropriate weights for instances of this kind. In the case of the modified Flower Snarks, the single added edge must be used in all Hamiltonian cycles (since the graph was non-Hamiltonian before its addition). This edge appears to always have the lowest resistance-distance out of all edges in the graph, which perhaps explains the strong performance. It appears that \lq\lq important" edges often have outstanding resistance-distance, although this can be in a minimal or maximal sense.

The {\em Minimal Reg} instances are a strong example of the kind of graphs which contain many Hamiltonian cycles, but contain some complex structure. Choosing to give each edge equal weight appears to have been a very poor strategy, but since each graph contains many Hamiltonian cycles, choosing random weights worked relatively well. Indeed, for these instances, random weights turned out to be superior to resistance distances, but not to the inverse of resistance distances.

Given the strong performance obtained by using the inverse of resistance-distance, and inspired by the field of electronics where the inverse of electrical resistance is called conductivity, we propose that this useful metric be given the name {\em conductivity distance}. Intuitively, it seems reasonable that conductivity distance would be a good choice of weights, since a large resistance distance implies many paths go between two vertices, and hence there are potentially more ways of finding a Hamiltonian cycle using the corresponding edge.

Although we have only used resistance distance and conductivity distance to solve instances of HCP in this note, the impact on solving time is so dramatic that it hints at the potential for future research into using these metrics for sparse TSP instances. Specifically, by augmenting the edge weights of a TSP instance with the resistance distance or conductivity distance (with intelligently chosen coefficients to prevent the true edge weights from being swamped), it might be possible to help coax TSP solvers such as Concorde to find good solutions more rapidly. This development could be widely useful, since even non-sparse TSP instances (such as those based on the Euclidean distance between each pair of vertices) are typically first converted to a sparse instance by a sparsification heuristic.

\end{document}